\documentclass[10pt]{amsart}
\usepackage{a4,tikz,enumerate,mathbbol,fancyhdr,colonequals,mathrsfs}
\numberwithin{equation}{section}

\theoremstyle{definition}

\numberwithin{theo}{section}
\author{Delio Mugnolo}
\keywords{Operator matrices, semigroups of
operators, wave equations with acoustic boundary conditions.}
\subjclass[2000]{47D05, 47H20, 35L20}

\title[Abstract Acoustic Boundary Conditions]{Abstract wave equations with acoustic boundary conditions}

\email{delio.mugnolo@uni-ulm.de}
\thanks{This paper has been written while I was a Ph.D. student at the University of T\"ubingen. I take this occasion to thank the Istituto Nazionale di Alta Matematica ``Francesco Severi'' for financial support as well as my supervisor, Rainer Nagel, for motivating discussions.}
\address{Abteilung Angewandte Analysis, Universit\"at Ulm, D-89081 Ulm, Germany}

\def\:{\thinspace:\thinspace}
\def\linie{\vrule height 14pt depth 5pt}
\def\back{\noalign{\vskip-3pt}}

\begin{document}

\begin{abstract}
We define an abstract setting to treat
wave equations equipped with time-dependent acoustic boundary
conditions on open domains of ${\bf R}^n$. We prove a
well-posedness result and develop a spectral theory which also allows
to prove a conjecture proposed in \cite{[GGG03]}. Concrete problems are also discussed.
\end{abstract}

\maketitle

\section{Introduction}

Wave equations equipped with homogeneous boundary conditions have been
studied for a long time. However, other kinds of boundary conditions
can also be considered, and for a number of concrete application it
seems that the right boundary conditions to impose are time-dependent,
cf. \cite{[Kr61]} and \cite{[BR74]}.

\bigskip
Certain investigations have in fact led theoreticyal physicists,
cf. \cite{[MI68]}, to investigate wave equations equipped with {\sl acoustic}
(or {\sl absorbing}) boundary conditions, which can be written in
the form
\begin{equation}\tag{ABC}
\left\{
\begin{array}{rcll}
\ddot{\phi}(t,x)&=& c^2 \Delta \phi(t,x),   &t\in {\bf R},\; x\in\Omega,\\
m(z)\ddot{\delta}(t,z)&=& -d(z)\dot{\delta}(t,z)-k(z){\delta}(t,z)-\rho(z)
\dot{\phi}(t,z), &t\in {\bf{R}},\;z\in\partial \Omega,\\
\dot{\delta}(t,z)&=& \frac{\partial\phi}{\partial\nu}(t,z),
&t\in {\bf R},\;z\in\partial \Omega.
\end{array}
\right.
\end{equation}
Here $\phi$ is the velocity potential of a fluid filling an open
 domain $\Omega\subset {\bf R}^n$, $n\geq1$; $\delta$ is the
normal displacement of the (sufficiently smooth) boundary $\partial \Omega$ of
 $\Omega$; $m$, $d$, and $k$ are the mass per unit area, the
 resistivity, and the spring constant of the boundary,
 respectively; finally, $\rho$ and $c$ are the unperturbed density of,
 and the constant speed of sound in the medium, respectively. It is
 reasonable to assume all these physical quantities to be modelled by
 essentially bounded functions, with $\rho,m$ real valued and  $\inf\limits_{z\in\partial \Omega} m(z)>0$. 

\bigskip
Quoting J.T. Beale and S.I. Rosencrans \cite{[BR74]} (who denote by $G$ our
domain $\Omega$), we point out that 'the physical model giving rise to
these conditions is that of a gas undergoing small irrotational
perturbations from rest in a domain $G$ with smooth compact boundary',
assuming that 'each point of the surface $\partial G$ acts like a
spring in response to the excess pressure in the gas, and that there
is no transverse tension between neighboring points of $\partial G$,
i.e., the ``springs'' are independent of each other'.

\bigskip
Operator matrices techniques have been used in this context already in
the 1970s, in a series of papers mainly by Beale. The well-posedness
of the initial value problem associated with $\rm(ABC)$ has
been announced in \cite{[BR74]}, and a detailed proof has been published
shortly afterwards (\cite[Thm.~2.1]{[Be76]}) under regularity assumptions on
the coefficients that are slightly more restrictive than ours.

Recently, acoustic boundary conditions have aroused interest
again. For example, C. Gal, G. Goldstein, and J.A. Goldstein
have compared them in \cite{[GGG03]} to more usual dynamical boundary
conditions for a wave equation, proving some spectral results and
proposing a conjecture that we are now able to prove.

It is remarkable that, back in the 1960s, the Russian school was
developing a spectral theory for extremely similar ``boundary-contact''
problems, whose most peculiar characteristic is that they possibly
feature in the boundary conditions differential operators of an order
that is higher than those acting in the interior, cf. \cite{[Kr61]}.

In the nice survey \cite{[Bl00]}, B. Belinsky has considered such
boundary-contact problems in order to describe some variations on an evocative
geophysical model. He has shown that the system of the generalized
eigenfunctions of an associated Sturm--Liouville problem
forms a Riesz basis on a suitable Sobolev space.

Let us also note that wave equations equipped with simpler boundary
conditions but much more complicated coupling relation have been
considered by G. Propst and J. Pr\"uss, who proved a well-posedness
result in \cite[Thm.~4.1]{[PP96]}.

\bigskip
Our purpose is to present a more general approach to such problems
that is based on results on operator matrices with non-diagonal
domain mainly obtained by K.-J. Engel (see, e.g., \cite{[En99]},
\cite{[CENN03]} and \cite{[KMN03]}). This reduces the need for formal
computations and allows more general cases, where the conditions of the
Lumer--Phillips theorem are harder to check.

Our paper is organized as follows. In Section~2 we introduce the
abstract setting we will exploit, and then show the well-posedness of
the abstract initial value problem associated with
(ABC). In Section~3 we sharpen some known results about Dirichlet
operators, and thereby investigate some spectral properties of a certain
operator matrix arising in our context. In particular, a conjecture
formulated in \cite{[GGG03]} is proven. In Section~4 we consider a
special case where the acoustic boundary conditions shrink to
dynamical boundary conditions of first order. This is of indpendent
interest, cf. \cite{CENP05}. Finally, motivated by a so-called Timoshenko
model discussed in \cite[\S~3]{[Bl00]}, we prove in Section~5 a well-posedness
result for second order problems with neutral acoustic boundary
conditions.

\section{General setting and well-posedness}

Inspired by the setting in \cite{[CENN03]}, we impose the following
throughout our paper.

\bigskip
\goodbreak
{\bf Assumptions~2.1.} 

(A$_1$) $X$, $Y$, and ${\partial X}$ are Banach spaces with
$Y\hookrightarrow X$.

(A$_2$) The operator $A:D(A)\to X$ is  linear with
$D(A)\subset Y$.

(A$_3$) The operator $R:D(A)\to {\partial X}$ is  linear and surjective.

(A$_4$) The operators $B_1,B_2$ are linear and bounded from $Y$
to ${\partial X}$.

(A$_5$) The operators $B_3,B_4$ are linear and bounded on ${\partial X}$.

(A$_6$) The operator ${A\choose R}:D(A)\subset Y\to X\times {\partial X}$
is closed.

(A$_7$) The restriction $A_0:=A_{\arrowvert\!\ker R}$ generates
a cosine operator function with associated phase space $Y\times X$.

\bigskip
Moreover, it will be convenient to define a new operator
$$L:=R+B_2,\qquad L:D(A)\to{\partial X}.$$
We will see that in some applications the operator $L$ is in some
sense ``more natural'' than $R$. E.g., when we discuss the motivating
equation $\rm(ABC)$, the operator $B_2$ will be the
trace operator and $L$ the normal derivative, while $R$ is a linear
combination of the two. This shows that the operator
$A_0=A_{\arrowvert\!\ker R}$ can be considered as an abstract version
of a operator equipped with Robin boundary conditions. (Recall that,
in the context of PDE's, Robin boundary conditions stand for boundary
conditions which are a linear combination of Dirichlet and Neumann
conditions, cf. \cite[\S~VII.3.2]{[DL90]}.) The main purpose of this paper is
to derive some properties of the wave equation with acoustic boundary
conditions from analogous properties of the wave equation with
homogeneous Robin (instead of Neumann, as in \cite{[Be76]}) boundary
conditions.

\bigskip
\goodbreak
{\bf Remark~2.2.}
Observe in particular that, due to the boundedness of $B_2$, the
condition (A$_6$) is satisfied if (and only if) also the operator
$${A\choose L}={A\choose R}+{0\choose B_2}:D(A)\subset Y\to X\times
{\partial X}$$
is closed.

\bigskip
Of concern in this paper are abstract second order initial-boundary
value problems equipped with (abstract) {\sl acoustic boundary
conditions} of the form
\begin{equation}\tag{AIBVP$_2$}
\left\{
\begin{array}{rcll}
 \ddot{u}(t)&=& Au(t), &t\in {\bf R},\\
 \ddot{x}(t)&=& B_1u(t)+B_2\dot{u}(t)+ B_3 x(t) + B_4\dot{x}(t),
 &t\in {\bf R},\\
  \dot{x}(t)&=& Lu(t), &t\in {\bf R},\\
  u(0)&=&f, \qquad\dot{u}(0)=g,&\\
  x(0)&=&h, \qquad\dot{x}(0)=j,&\\
\end{array}
\right.
\end{equation}
on $X$ and ${\partial X}$, where the operators $A,B_1,B_2,B_3,B_4$, and
 $R=L-B_2$ satisfy the Assumptions~2.1.

\bigskip
We want to recast the second order problem ${\rm(AIBVP_2)}$ as a
first order initial-boundary value problem. Such problems
 have been thoroughly discussed in \cite{[KMN03]}. This approach is mostly
 based on the notion of, and on some results on so-called {\sl
one-sided coupled operator matrices}, cf.\cite{[En99]} and \cite{[CENN03]},
 and exploits semigroup theory as an essential instrument. In fact, the
results in this section are strictly related to properties of the
operator matrix with non-diagonal domain ${\mathcal{A}}$ as defined in (3.1).

Thus, we re-write $\rm(AIBVP_2)$ as a
first order abstract initial-boundary value problem
\begin{equation}\tag{$\mathbb{AIBVP}$}
\left\{
\begin{array}{rcll}
 \dot{\mathbb{u}}(t)&=& \mathbb{A}\mathbb{u}(t), &t\in {\bf R},\\
 \dot{\mathbb{x}}(t)&=&
 \mathbb{B}\mathbb{u}(t)+\tilde{\mathbb{B}}\mathbb{x}(t), &t\in {\bf R},\\ 
 \mathbb{x}(t)&=&\mathbb{R}\mathbb{u}(t), &t\in {\bf R},\\
  \mathbb{u}(0)&=&\mathbb{u}_0, &\\
  \mathbb{x}(0)&=&\mathbb{x}_0, &
\end{array}
\right.
\end{equation}
on the Banach spaces
$$\mathbb{X} :=Y\times X\times {\partial X}\qquad\hbox{and}\qquad\partial
\mathbb{X} :={\partial X}.\leqno(2.1)$$
The operator $\mathbb{A} $ on $\mathbb{X} $ is given by
$${{\mathbb{A} }}:=\begin{pmatrix}
0 & I_Y  & 0\\
A & 0    & 0\\
L & 0    & 0
\end{pmatrix},\qquad
D(\mathbb{A} ):=D(A)\times Y\times {\partial X}.\leqno(2.2)$$
Further, $\mathbb{R} $ and $\mathbb{B} $ are the operators
$$\mathbb{R} :=\begin{pmatrix}
R & 0 & 0
\end{pmatrix},
\qquad D(\mathbb{R} ):=D(\mathbb{A} ),\leqno(2.3)$$
and
$$\mathbb{B} :=\begin{pmatrix}
B_1+B_4B_2 & 0 & B_3
\end{pmatrix},
\qquad D(\mathbb{B} ):=\mathbb{X} ,\leqno(2.4)$$
respectively, both from $\mathbb{X} $ to $\partial\mathbb{X}$. Moreover,
$\tilde{\mathbb{B} }$ is the operator
$$\tilde{\mathbb{B} }:=B_4,\qquad D(\tilde{\mathbb{B} }):=\partial \mathbb{X}
,\leqno(2.5)$$ 
on $\partial\mathbb{X} $.
Finally, we set the initial data
$$\mathbb{u}_0:=\begin{pmatrix}
f\\ g\\ h
\end{pmatrix}
\qquad\hbox{and}\qquad
\mathbb{x}_0:=j-B_2f.\leqno(2.6)$$

\bigskip
\goodbreak
{\bf Lemma~2.3}
{\sl 
The following assertions hold.

- The restriction $\mathbb{A}_{0}:=\mathbb{A}_{\arrowvert\!\ker \mathbb{R} }$ generates
a strongly continuous group on $\mathbb{X} $.

- The operator $\mathbb{R} $ is surjective.

- The operator $\mathbb{B} $ is bounded from $\mathbb{X} $ to $\partial \mathbb{X}$.

- The operator $\tilde{\mathbb{B} }$ is bounded on $\partial \mathbb{X} $.

- The operator ${\mathbb{A} \choose \mathbb{R} }:D(\mathbb{A} )\subset
\mathbb{X} \to \partial \mathbb{X} $ is closed.}

\begin{proof}
Observe first that $\ker \mathbb{R} =\left\{u\in D(A)\: Lu=B_2u\right\}
\times Y\times {\partial X}$, thus the operator $\mathbb{A}_0$ takes the form
$$\mathbb{A}_0=\begin{pmatrix}
\lineskip=0pt
0 & I_Y &\linie & 0\\ \back
A_0 & 0 &\linie & 0\\
\noalign{\hrule}
B_2 & 0 &\linie & 0
\end{pmatrix}.\leqno(2.7)$$
Observe that the perturbation $\begin{pmatrix}B_2 & 0\end{pmatrix}$ is bounded from
$Y\times X$ to ${\partial X}$, and the only non-zero diagonal block of $\mathbb{A}_0$
generates by \cite[Thm.~3.14.11]{[ABHN01]} a strongly continuous group on
$Y\times X$.
Therefore, $\mathbb{A}_0$ generates a strongly continuous group on $\mathbb{X} $, and
                (i) is proven. 
The remaining claims follow by Assumptions~(A$_3$)--(A$_6$).
\end{proof}

\bigskip
Therefore, by \cite[Prop.~4.1]{[KMN03]}, the following result is immediate.

\bigskip
\goodbreak
{\bf Proposition~2.4.}
{\sl The abstract initial-boundary value problem
$(\mathbb{AIBVP})$ is well-posed in the sense of \cite[\S~2]{[KMN03]}.}

\bigskip
We now come back to the discussion of the original second order
abstract initial-boundary value problem $\rm{(AIBVP_2)}$.

\bigskip
\goodbreak
{\bf Definition~2.5.}
A function $u:{\bf R}\to X$ is a
classical solution to $\rm(AIBVP_2)$ on $(Y,X,{\partial X})$ if

- $u\in C^2({\bf R}, X)\cap C^1({\bf R}, Y) $ and $Lu\in
C^1({\bf R}, {\partial X})$,

- $u(t)\in D(A)$ for all $t\in {\bf R}$,

- $u$ satisfies $\rm(AIBVP_2)$ pointwise.

\bigskip
We will identify solutions to $\rm(AIBVP_2)$ on
$(Y,X,{\partial X})$ and solutions to $(\mathbb{AIBVP})$ by letting
$$\mathbb{u} (t)\equiv 
\begin{pmatrix}
u(t)\\ \dot{u}(t)\\ h+\int_0^t Lu(s)\;ds
\end{pmatrix}
\qquad \hbox{and}\qquad \mathbb{x} (t)\equiv Ru(t),\qquad
t\in {\bf R}.$$
This is justified by the following.

\bigskip
\goodbreak
{\bf Lemma~2.6.} 
{\sl The problems $\rm(AIBVP_2)$ and $(\mathbb{AIBVP})$
are equivalent.}

\begin{proof} 
Let
$$\mathbb{u} =\begin{pmatrix}
u\\ v\\ x
\end{pmatrix}
\in C^1({\bf R},\mathbb{X} )$$
be a classical solution to $(\mathbb{AIBVP})$ and let $\mathbb{R} \mathbb{u} =Ru=y$. Thus,
there holds
\begin{equation*}
\left\{
\begin{array}{rcll}
 \dot{u}(t)&=&v(t), &t\in {\bf R},\\
 \dot{v}(t)&=&Au(t), &t\in {\bf R},\\
 \dot{x}(t)&=&Lu(t), &t\in {\bf R},\\
 \dot{y}(t)&=& (B_1+B_4 B_2) u(t)
 +B_3 x(t) + B_4 y(t), &t\in {\bf R},\\
 y(t)&=&(L-B_2)u(t), t\in {\bf R},\\
 u(0)&=&f, \qquad v(0)=g,&\\
 x(0)&=&h, \qquad y(0)=j-B_2f,&
\end{array}
\right.
\end{equation*}
with $v(t)\in Y$, $t\in {\bf R}$, or, equivalently,
\begin{equation*}
\left\{
\begin{array}{rcll}
 \ddot{u}(t)&=&Au(t), &t\in {\bf R},\\
 \dot{x}(t)&=&Lu(t), &t\in {\bf R},\\
  \ddot{x}(t)&=& B_1 u(t) +B_2\dot{u}(t)+ B_3 x(t) +
 B_4\dot{x}(t), &t\in {\bf R},\\
 u(0)&=&f, \qquad \dot{u}(0)=g,&\\
 x(0)&=&h, \qquad \dot{x}(0)=j,&
\end{array}
\right.
\end{equation*}
with $\dot{u}(t)\in Y$, $t\in {\bf R}$. To justify this step observe
that, by assumption,
 $u\in C^1({\bf R},Y)$. Therefore we see that
$$B_2 \dot{u}(t)=B_2\left(Y\!\!-\!\!\lim_{h\to 0}{u(t+h)-u(t)\over
h}\right)={\partial X}\!\!-\!\!\lim_{h\to 0} B_2\left({u(t+h)-u(t)\over
h}\right)={d\over dt}B_2 u(t),$$
where we have used the assumption $B_2\in{\mathcal L}(Y,{\partial X})$. Note
that this argument does not hold for $L$.

This shows that $u$ is actually a classical solution to $\rm(AIBVP_2)$
on $(Y,X,{\partial X})$. The converse implication follows likewise, and the
claim is proven.
\end{proof}

\bigskip
Once we have shown the well-posedness of $\rm(AIBVP_2)$, we can look
back at the original initial value problem associated with the wave
equation $\rm(ABC)$ introduced in Section~1. Thus, we obtain the
following.

\bigskip
{\bf Theorem~2.7.}
{\sl The initial value problem associated
with the wave equation with acoustic boundary conditions $\rm(ABC)$
on an open domain $\Omega\subset {\bf R}^n$, $n\geq1$, with
smooth boundary $\partial \Omega$ is well-posed. In particular, for all
initial data 
$$\phi(0,\cdot)\in H^2(\Omega),\qquad
\dot{\phi}(0,\cdot)\in H^1(\Omega),\qquad \delta(0,\cdot)\in
L^2(\partial \Omega),$$
$$\hbox{and}\qquad\dot{\delta}(0,\cdot)\in L^2(\partial \Omega)
\qquad\hbox{such that}\qquad
{\partial \phi\over \partial \nu}(0,\cdot)=\dot{\delta}(0,\cdot)$$
there exists a classical solution on
$(H^1(\Omega),L^2(\Omega),L^2(\partial \Omega))$ continuously depending on
them.}

\begin{proof}
Take first 
$$X:=L^2(\Omega),\qquad Y:=H^1(\Omega),\qquad
{\partial X}:=L^2(\partial \Omega).$$ 
We set
$$A:=c^2\Delta,\qquad D(A):=\left\{u\in H^{3\over 2}(\Omega)\: 
\Delta u\in L^2(\Omega)\right\},$$
$$(Rf)(z)={\partial f\over\partial \nu}(z)+{\rho(z)\over
m(z)} f(z),\qquad f\in D(R)=D(A),\; z\in\partial \Omega,$$
$$B_1=0,\qquad (B_2 f)(z):=-{\rho(z)\over m(z)} f(z),\qquad 
f\in H^1(\Omega),\; z\in\partial \Omega,$$
$$(B_3 g)(z):=-{k(z)\over m(z)}g(z),\qquad (B_4 g)(z):=-{d(z)\over
  m(z)}g(z),\qquad g\in L^2(\partial \Omega),\; z\in\partial \Omega.$$
By Proposition~2.4 and Lemma~2.6, it suffices to check that the
Assumptions~(A$_1$)--(A$_7$) are satisfied in the above
setting.

The Assumptions~(A$_1$) and (A$_2$) are clearly satisfied. To
check the Assumption~(A$_3$), we apply \cite[Vol.~I, Thm.~2.7.4]{[LM72]} and
obtain that for all $g\in L^2(\partial \Omega)$ there exists a $u\in H^{3\over
2}(\Omega)$ such that $\Delta u=0$ and ${\partial u\over \partial
\nu}+{\rho\over m}u_{\arrowvert\partial \Omega}=g$. The Assumption~(A$_4$)
holds because the trace operator is bounded from $H^1(\Omega)$ to
$L^2(\partial \Omega)$ 
and because ${\rho\over m}\in L^\infty(\partial \Omega)$, while the
Assumption~(A$_5$) is satisfied since ${d\over m},{k\over m}\in
L^\infty(\partial \Omega)$.

The Assumption~(A$_6$) is satisfied because the closedness of
${A\choose L}$ holds by interior estimates for elliptic operators,
(a short proof of this can be found in \cite[\S~3]{[CENN03]}),
and $B_2\in{\mathcal L}(Y,{\partial X})$, cf. Remark~2.2.

To check Assumption~(A$_7$), observe that the operator
$A_0=A_{\arrowvert\!\ker R}$ is in fact (up to the constant $c^2$) the
Laplacian with Robin boundary conditions, that is,
$$A_0 u = c^2 \Delta u,\qquad D(A_0)=\left\{u\in H^2(\Omega)\:
{\partial   u\over\partial\nu}+{\rho\over m}
u_{\arrowvert\partial \Omega}=0\right\}.$$ 
This operators is self-adjoint and dissipative up to a scalar perturbation. By the results of~\cite[\S~7.1]{ABHN01}, it  generates a cosine operator function with associated phase space $H^1(\Omega)\times L^2(\Omega)= Y\times X$.
\end{proof}

\bigskip
Similar well-posedness results for the initial value problem
associated with (ABC) have already been obtained in \cite{[Be76]}.
However, the coefficients $\rho,d,k,m$ are therein assumed to satisfy
more restrictive assumptions, as in particular they are supposed to be
real and positive ($m$ strictly positive, $\rho$ constant) continuous
functions on $\partial \Omega$.

Our Assumptions~2.1 are satisfied by a variety of other operators and
spaces. We discuss a biharmonic wave equation with acoustic-type boundary
conditions.

\bigskip
{\bf Example~2.8.} Let $p,q,r,s\in L^\infty(\partial \Omega)$. Then the initial value problem associated with
\begin{equation*}
\left\{
\begin{array}{rcll}
 \ddot{\phi}(t,x)&=& -\Delta^2 \phi(t,x), &t\in {\bf R},\;
 x\in \Omega,\\
\ddot{\delta}(t,z)&=&  p(z){\delta}(t,z)+q(z)\dot{\delta}(t,z)+r(z)
{\partial \phi\over\partial \nu}(t,z)+s(z)
\frac{\partial \dot{\phi}}{\partial \nu}(t,z), & t\in {\bf
R},\;z\in\partial \Omega,\\
\dot{\delta}(t,z)&=&{\Delta\phi}(t,z), & t\in {\bf R},\;z\in\partial \Omega,\\
\phi(t,z)&=&0, &t\in {\bf R},\;z\in\partial \Omega,
\end{array}
\right.
\end{equation*}
is well-posed. In particular, for all initial data
$$\phi(0,\cdot)\in H^4(\Omega)\cap H^1_0(\Omega),\qquad
\dot{\phi}(0,\cdot)\in H^2(\Omega)\cap H^1_0(\Omega),\qquad
\delta(0,\cdot)\in L^2(\partial \Omega),$$
$$\hbox{and}\qquad\dot{\delta}(0,\cdot)\in L^2(\partial \Omega)\qquad\hbox{such
that}\qquad{\partial^2
\phi\over\partial\nu^2}(0,\cdot)=\dot{\delta}(0,\cdot)$$ 
there exists a classical solution continuously depending on them.

Take 
$$X:=L^2(\Omega),\qquad
Y:=H^2(\Omega)\cap H^1_0(\Omega),\qquad
{\partial X}:=L^2(\partial \Omega),$$
 and consider the operators
$$A:=-\Delta^2 ,\qquad 
D(A):=\{u\in H^{5\over 2}(\Omega)\cap H^1_0(\Omega): \Delta^2 u\in
L^2(\Omega)\},$$
$$Ru:=\left(\Delta u\right)\arrowvert_{\partial\Omega} -s{\partial u \over \partial \nu},\qquad \hbox{for all}\quad u\in D(R):=D(A),$$
$$B_1:=r{\partial \over\partial \nu},\qquad B_2:=s{\partial 
\over\partial \nu},\qquad D(B_1):=D(B_2):=Y,$$
$$B_3x:=px,\qquad B_4:=qx,\qquad \hbox{for all}\quad x\in {\partial X}.$$
We are only going to prove that $A_0$, i.e., the restriction of
$-\Delta^2$ to
$$D(A_0):=\ker R=\left\{u\in H^4(\Omega)\cap H^1_0(\Omega)\: 
\left(\Delta u\right)\arrowvert_{\partial\Omega} =
s{\partial u\over \partial\nu}\right\},$$ 
generates a cosine operator function with associated phase space
$\big(H^2(\Omega)\cap H^1_0(\Omega)\big)\times L^2(\Omega)=Y\times
X$, the remaining Assumptions~2.1 being satisfied trivially.

Take $u,v\in D(A_0)$ and observe that applying the
Gauss-Green formulas twice yields
$$\big<A_0 u,v\big>_X=-\int_\Omega \Delta^2 u \cdot\overline{v}\;
dx=-\int_\Omega \Delta u\cdot\overline{\Delta v}\; dx +\int_{\partial \Omega} s\;
{\partial u\over \partial \nu} \cdot\overline{\partial v\over \partial
\nu}\; d\sigma.$$
It is immediate that $A_0$ is self-adjoint and dissipative up to a scalar perturbation, hence by the results of~\cite[\S~7.1]{ABHN01}, the generator of a cosine
operator function with associated phase space $V\times X$, for some
Banach space $V$. We claim that $V=Y=H^2(\Omega)\cap H^1_0(\Omega)$,
thus that the associated phase space is actually $Y\times X$.

Integrating by parts one 
sees that $0$ is not an eigenvalue of $A_0-\omega$, for $\omega>0$ large enough, hence $-A_0+\omega$ is a strictly positive self-adjoint operator. The domain of its square
root coincides by \cite[Thm.~VI.2.23]{[Ka95]} with the form domain
of $A_0$. Moreover, one can directly check that the form domain of
$A_0$ is $H^2(\Omega)\cap H^1_0(\Omega)$. Therefore, by
\cite[Prop.~VI.3.14]{[EN00]} we deduce that the associated
phase space of $A_0$ is $V\times X=[D(-A_0+\omega)^{1\over 2}]\times X=
\big(H^2(\Omega)\cap H^1_0(\Omega)\big)\times  L^2(\Omega)$ as
claimed.

\bigskip
{\bf Remark~2.9.} Among further operators and spaces fitting into 
our abstract framework we list the following. In both
cases, the operator $B_2$ is defined as in the proof of Theorem~2.7.

a) $X:=L^2(\Omega)$, $\;Y:=H^1(\Omega)$, $\;{\partial X}:=L^2(\partial
\Omega)$, 

\qquad $Au(x):=\nabla( a(x)\nabla u(x))$, $\; x\in\Omega$, with the
function $a\ge0$ sufficiently regular on $\overline\Omega$,

\qquad $Lu(z)=<a(z)\nabla u(z),\nu(z)>$, $\; z\in\partial \Omega$,

\qquad for $u\in D(A):=\{H^{3\over 2}(\Omega):Au\in L^2(\Omega)\}$.

b) $X:=L^2(\Omega)$, $\;Y=H^2(\Omega)$, $\;{\partial X}:=L^2(\partial \Omega)$,

\qquad $Au:=-\Delta^2 u$,

\qquad $Lu:=-{\partial \Delta u\over\partial \nu}$,

\qquad for $u\in D(A):=\left\{H^{7\over 2}(\Omega): \Delta^2 u\in
L^2(\Omega),\; \left( \Delta u\right)\arrowvert_{\partial\Omega}=0\right\}$.

\bigskip
\noindent
In either case, $A_0$ is self-adjoint and dissipative up to a scalar perturbation. This ensures that $A_0$ generates a cosine operator function.

\bigskip
Before concluding this section, let us emphasize that 
our proof of the well-posedness of $\rm(AIBVP_2)$ actually 
relies on the reformulation of $(\mathbb{AIBVP})$
as a first order abstract Cauchy problem
\begin{equation}\tag{$\mathcal{ACP}$}
\left\{
\begin{array}{rcll}
\dot{\mathcal{U}}(t)&=&{\mathcal{A}}\;{\mathcal{U}}(t), &t\in {\bf R},\\
{\mathcal{U}}(0)&=&{\mathcal{U}}_0,
\end{array}
\right.
\end{equation}
on ${\mathcal{X}}:=\mathbb{X} \times \partial \mathbb{X} $, where
$${\mathcal{A}}:=\begin{pmatrix}
\mathbb{A}  & 0 \\ 
\mathbb{B}  & \tilde{\mathbb{B} }
\end{pmatrix},\qquad
D({\mathcal{A}}):=\left\{{\mathbb{u} \choose \mathbb{x} }\in D(\mathbb{A}
)\times \partial \mathbb{X} \: 
\mathbb{R} \mathbb{u} =\mathbb{x} \right\},\leqno(2.8)$$
is an operator on ${\mathcal{X}}$ and
$${\mathcal{U}}(t):={\mathbb{u} (t)\choose \mathbb{R} \mathbb{u} (t)},\quad
t\in {\bf R},\qquad 
{\mathcal{U}}_0:={\mathbb{u}_0\choose \mathbb{x}_0}.$$
In fact, it has beens shown in \cite[\S~2]{[KMN03]} that a
function $\mathbb{u} :{\bf R}\to \mathbb{X} $ is a classical solution to
$(\mathbb{AIBVP})$ if it is the first coordinate of a classical solution
${\mathcal{U}}:{\bf R}\to {\mathcal{X}}$ to $(\mathcal{ACP})$, and in fact the
proof of Proposition~2.4 relies on the fact that, by Lemma~2.3, 
${\mathcal{A}}$ generates a group on $\mathcal{X}$,
cf. \cite[Prop.~4.1]{[KMN03]}.

\bigskip 
{\bf Remark~2.10.} Observe that, in order to show the well-posedness
of (a problem equivalent to) $(\mathcal{ACP})$ on an open bounded domain $\Omega$ of ${\bf R}^3$, Beale introduced the product space 
$$(H^1(\Omega;\rho)/{\bf C})\times
L^2(\Omega;\frac{\rho}{c^2})\times L^2(\partial \Omega;k)\times
L^2(\partial \Omega;m).$$ 
Such a space looks somehow artificial, due to the quotient appearing in the first coordinate and to the weights of the remaining $L^2$-spaces.
He then showed that a certain operator matrix (different from our
${\mathcal{A}}$) verifies the conditions of the Lumer--Phillips theorem --
that is, the energy of the solutions to $\rm(ABC)$ is nonincreasing for time
$t\ge0$. Moreover, if the parameter $d\equiv 0$, then also the
conditions of Stone's theorem are satisfied -- that is, the energy
is constant for $t\in {\bf R}$. Also, Beale showed that his operator matrix does not have compact resolvent and computed its essential spectrum, but his techniques can hardly be applied to problems on domains of ${\bf R}^n$, $n\not=3$. 

Since Beale's paper, the theory of asymptotics for (semi)groups has been widely developped. In particular, it is now known that every bounded strongly continuous group (resp., semigroup) whose generator has compact resolvent is almost periodic (resp., asymptotically almost periodic), cf.~\cite[Chapt.~5]{[ABHN01]}. (More generally, every bounded strongly continuous semigroup with only countably many spectral values on $i\bf R$ is asymptotically almost periodic. It seems therefore worthwhile to develop a complete spectral theory for the problem (AIBVP$_2$).) This is the main aim of Section~3.

The main drawback of our own approach is that we fail to produce an
energy space for the motivating equation $\rm(ABC)$ on a bounded domain
$\Omega\subset {\bf R}^n$, i.e., the group generated by ${\mathcal{A}}$ is not 
contractive, as it can be seen already in the case of $n=1$.

However, our approach has other advantages. In particular, the above
operator matrix ${\mathcal{A}}$ can be written as
$${\mathcal{A}}:={\mathcal{A}}_1+{\mathcal{A}}_2:= 
\begin{pmatrix}
\mathbb{A}  & 0\\ 
0 & 0
\end{pmatrix}+
\begin{pmatrix}
0 & 0\\
\mathbb{B}  & \tilde{\mathbb{B} }
\end{pmatrix},\leqno(2.9)$$ 
where ${\mathcal{A}}_1$ is the generator of a strongly continuous group on
${\mathcal{X}}$  
and ${\mathcal{A}}_2$ is a bounded operator, which is compact if and only if 
$\rm dim\;{\partial X}<\infty$. 
A decomposition of this type cannot be performed on the operator 
matrix considered by Beale, and has some interesting consequences:
Observe that, in the context of $\rm(ABC)$, the group generated by 
${\mathcal{A}}_1$ governs the inital value problem associated with the wave 
equation with inhomogeneous (static) Robin boundary conditions
\begin{equation}\tag{iRBC}
\left\{
\begin{array}{ll}
\ddot{\psi}(t,x)= c^2 \Delta \psi(t,x), &t\in {\bf R},\;
 x\in \Omega,\\
{\partial\psi \over\partial\nu}(t,z)+{\rho(z)\over m(z)}\psi(t,z)=
\gamma(z), 
&t\in {\bf R},\;z\in\partial \Omega,\\
\end{array}
\right.
\end{equation}
where $\gamma(z):={\partial\psi\over\partial \nu}(0,z)+{\rho(z)\over 
 m(z)}\psi(0,z)$.

The solution $\psi$ to $\rm(iRBC)$ is given by the group generated by
 ${\mathcal{A}}_1$ and can be explicitly written down, cf. 
\cite[Thm.~3.6]{[KMN03]}. Thus, due to (2.9), the solution $\phi$ to
 $\rm(ABC)$ can be obtained by the Dyson--Phillips series,
 cf. \cite[Thm.~III.1.10]{[EN00]}. Further, by \cite[Cor.~III.1.11]{[EN00]} we
 deduce that  
$$\| \phi(t,\cdot)-\psi(t,\cdot)\|_{L^2(\Omega)}\leq tM,$$
for $t\in[0,1]$ and some constant $M$.

\section{Regularity and spectral theory}

Due to the important role played by the operator matrix ${\mathcal{A}}$
defined in (2.8), we are interested in developing a spectral theory for it. To
this purpose we are still imposing the Assumptions~2.1.

As a first step, we recall the notion of Dirichlet operators. To
obtain an optimal boundedness result, we need the following.

\bigskip
{\bf Lemma~3.1.}
{\sl Let $Z$ be a Banach space such that
$Z\hookrightarrow Y$, and consider the operator matrix
$$\begin{pmatrix}
A & 0\\ 
R & 0
\end{pmatrix}
:D(A)\times {\partial X}\to X\times {\partial X}$$
on $X\times {\partial X}$.
Then its part in $Z\times {\partial X}$ is closed.}

\begin{proof}
We can consider the part $A_\arrowvert$ of $A$ in $Z$ and let
${u_n\choose x_n}_{n\in {\bf N}}\!\!\subset D(A_\arrowvert)\times
{\partial X}$ such that
$$ {u_n\choose x_n}\to {u\choose x}\qquad\hbox{in}\; Z\times{\partial X}$$
and
$$\begin{pmatrix}
A_\arrowvert & 0\\ 
R & 0
\end{pmatrix}
\begin{pmatrix}
u_n \\ x_n
\end{pmatrix}=
\begin{pmatrix}
Au_n \\ Ru_n
\end{pmatrix}
\to 
\begin{pmatrix}
w \\ y
\end{pmatrix}\qquad\hbox{in}\; Z\times{\partial X}.$$
Since $Z\hookrightarrow Y$, we can apply the closedness of $A\choose R$
and conclude that $u\in D(A)$, $Au=w$, and $Ru=y$. This completes the
proof.
\end{proof}

\bigskip
We are now in the position to show the following refinement of
\cite[Lemma~2.2]{[CENN03]}.

\bigskip
{\bf Lemma~3.2.} {\sl If $\lambda\in\rho(A_0)$, then the restriction
 $R\big\arrowvert_{\ker(\lambda-{A})}$ has an inverse
$$D^{A,R}_\lambda :\partial X\to\ker(\lambda-{A}),$$
called Dirichlet operator associated with $A$ and $R$.
Moreover, $D_\lambda^{A,R}$ is bounded from ${\partial X}$ to $Z$ for every
Banach space $Z$ satisfying $D(A^\infty)\subset Z\hookrightarrow
Y$.}

\begin{proof} The existence of the Dirichlet operator $D^{A,R}_\lambda$
follows from \cite[Lemma~2.2]{[CENN03]}, due to the
Assumptions~(A$_3$),(A$_7$).

Observe now that $\ker(\lambda-A)\subset D(A^\infty)$. Therefore the
 boundedness of $D_\lambda^{A,R}$ from ${\partial X}$ to some Banach space $Z$
 containing $D(A^\infty)$ is  equivalent to the closedness of the
 operator $R\big\arrowvert_{\ker(\lambda-{A})}:\ker(\lambda-A)\subset
 Z\to {\partial X}$.

To show that
$R\big\arrowvert_{\ker(\lambda-{A})}$ is actually closed, take
$(u_n)_{n\in {\bf N}}\subset\ker(\lambda-A)$ such that $u_n{\buildrel
Z\over \to}u$ and $Ru_n{\buildrel {\partial X}\over \to}x$. It follows that
$Au_n=\lambda u_n{\buildrel Z\over \to}\lambda u$, that is
$$\begin{pmatrix}
A & 0\\ 
R & 0
\end{pmatrix}_{\big\arrowvert}{u_n \choose
0}\to {\lambda u\choose x}\qquad\hbox{in}\; Z\times{\partial X}.$$
By Lemma~3.1 we conclude that $u\in D(A)$ and that $Au=\lambda u$,
$Ru=x$.
\end{proof}

\bigskip 
{\bf Remarks~3.3.} a) The above Dirichlet operators $D_\lambda^{A,R}$
can also be interpreted as follows: Lemma~3.2 says that the abstract
eigenvalue problem
\begin{equation*}
\left\{
\begin{array}{rcll}
Au&=& \lambda u &\hbox{in}\; X,\\
 Ru&=& x &\hbox{in}\; {\partial X},
\end{array}
\right.
\end{equation*}
has a unique solution $D_\lambda^{A,R}x$ for all $x\in{\partial X}$, and that 
the dependence on $x$ is continuous.

b) If ${\partial X}$ is finite dimensional, or else if a Banach space $Z$
as in the statement of Lemma~3.2 can be chosen to be {\sl compactly}
embedded in $Y$, then we obtain that
the Dirichlet operators are compact from ${\partial X}$ to $Y$.

c) We finally observe that, by definition,
$$LD_\lambda^{A,R}=I_{\partial X}+B_2 D_\lambda^{A,R}\qquad\hbox{for
all}\;\lambda\in\rho(A_0),\leqno(3.1)$$
and therefore $LD_\lambda^{A,R}$ is a bounded operator on ${\partial X}$.

\bigskip
We recall that the resolvent set of the operator matrix
$$\begin{pmatrix}
0 & I_Y \\ 
A_0 & 0
\end{pmatrix}\qquad\hbox{with domain}\qquad
D(A_0)\times Y\leqno(3.2)$$
on the space $Y\times X$ is given by 
$\{\lambda\in {\bf C}\:\lambda^2\in\rho(A_0)\}$. Accordingly, we 
obtain the following.

\bigskip
{\bf Lemma~3.4.} {\sl The resolvent set of $\mathbb{A}_0$ is given by
$$\rho(\mathbb{A}_0)=\left\{\lambda\in
{\bf C}\:\lambda\not=0,\;\lambda^2\in\rho(A_0)\right\}.$$
For $\lambda\in\rho(\mathbb{A}_0)$ there holds
$$R(\lambda,\mathbb{A}_0)=
\begin{pmatrix}
\lambda R(\lambda^2,A_0) & R(\lambda^2,A_0) & 0\\
A_0 R(\lambda^2,A_0) & \lambda R(\lambda^2,A_0) & 0\\
-B_2R(\lambda^2,A_0) & -{1\over \lambda}B_2R(\lambda^2,A_0) &
{1\over \lambda}I_{\partial X}
\end{pmatrix}.$$}

\begin{proof} The resolvent operator of the operator matrix introduced
in (3.2) is given by \cite[(3.107)]{[ABHN01]}. Then the claimed formula 
can be checked directly.
\end{proof}

\bigskip
{\bf Lemma~3.5.} {\sl For the operator $(\mathbb{A} ,D(\mathbb{A} ))$ defined in (2.2) we
obtain
$$D(\mathbb{A}^{2k-1})=D(A^k)\times D((A^{k-1})_{\arrowvert Y})
\times {\partial X}\qquad\hbox{and}$$
$$D(\mathbb{A} ^{2k})=D((A^k)_{\arrowvert Y})\times D(A^{k})\times
{\partial X}\qquad\hbox{for all}\; k\in {\bf N}.$$
In particular, $D(\mathbb{A} ^\infty)=D(A^\infty)\times D(A^\infty)\times
{\partial X}$.}

\begin{proof} The claim follows by induction on $n$,
using the fact that
$$D((A^k)_{\arrowvert Y})=\{u\in D(A)\: Au\in
D((A^{k-1})_{\arrowvert Y})\}$$
and recalling that $D(A)\subset Y$.
\end{proof}

\bigskip
This allows to obtain the following regularity result. 

\bigskip
{\bf Corollary~3.6.} {\sl Assume that the initial data $f,g$ are in 
$$\mathcal{D}^\infty_0:=\bigcap_{k=0}^\infty\{w\in D(A^k)\: RA^k w =LA^k w
=0\}.\leqno(3.3)$$   
If further $h=j=0$, then the solution $u=u(t)$ to $\rm(AIBVP_2)$ is in
$D(A^\infty)$ for all $t\in {\bf R}$.}

\begin{proof} Taking into account Lemma~3.5, the inclusion
  $(\mathcal{D}^\infty_0)^2 \times\{0\}^2\;\subset\; D({\mathcal{A}}^\infty)$
  can be proven by induction. Then, one only needs to recall that the group
  generated by ${\mathcal{A}}$ maps $D({\mathcal{A}}^n)$ in itself for all
  $n\in{\bf N}$, and to observe that $D({\mathcal{A}}^\infty)\subset
  D(A^\infty)^2\;\times\; {\partial X}^2$.
\end{proof}

\bigskip
{\bf Example~3.7.} Consider the framework introduced in the
proof of Theorem~2.7 to treat the motivating equation $\rm(ABC)$. 
Corollary~3.6 yields a regularity result that is similar to
\cite[Thm.~2.2]{[Be76]}. In fact, $\mathcal{D}_0^\infty$ defined in (3.3) contains
the set of all functions of class $C^\infty(\Omega)$ that vanish in a
suitable neighborhood of $\partial \Omega$. On the other hand,
$D(A^\infty)\subset H^{3\over 2}(\Omega)\cap C^\infty(\Omega)$.

\bigskip
Observe now that, by Lemma~2.3, the operators $\mathbb{A} $ and $\mathbb{R} $
satisfy the Assumptions~2.1. Therefore by Lemma~3.2 and for
$\lambda\in\rho(\mathbb{A}_0)$, one obtains the existence of the Dirichlet
operator $D_\lambda^{\mathbb{A} ,\mathbb{R} }$ associated with $\mathbb{A} $
and $\mathbb{R} $. More precisely, the following representation holds.

\bigskip
{\bf Lemma~3.8.} {\sl Let $\lambda\in\rho(\mathbb{A}_0)$. Then the Dirichlet
operator $D_\lambda^{\mathbb{A} ,\mathbb{R} }$ exists and is represented by
$$D_\lambda^{\mathbb{A} ,\mathbb{R} }=
\begin{pmatrix}
D_{\lambda^2}^{A,R}\\
\lambda D_{\lambda^2}^{A,R}\\
{1\over \lambda} L D_{\lambda^2}^{A,R}
\end{pmatrix}.\leqno(3.4)$$
Moreover, $D_\lambda^{\mathbb{A} ,\mathbb{R} }$ is bounded from $\partial
\mathbb{X} $ to $W\times Z\times {\partial X}$ for every two Banach spaces
$W,Z$ such that $D(A^\infty)\subset W\hookrightarrow Y$ and
$D(A^\infty)\subset Z\hookrightarrow X$.}

\begin{proof} To obtain the claimed representation, take
$\mathbb{x} :=y\in{\partial X}=\partial \mathbb{X} $. By definition the Dirichlet operator
$D_\lambda^{\mathbb{A} ,\mathbb{R} }$ maps $\mathbb{X} $ into the unique vector
\begin{equation*}
\mathbb{u} :=
\begin{pmatrix}
u\\ v\\ x
\end{pmatrix}
\in \mathbb{X} \qquad\hbox{such that}\qquad
\left\{
\begin{array}{rcl}
\mathbb{A} \mathbb{u} &=&\lambda \mathbb{u} ,\\
 \mathbb{R} \mathbb{u} &=&\mathbb{x} ,
\end{array}
\right.
\qquad\hbox{or rather}\qquad
\left\{
\begin{array}{rcl}
v&=&\lambda u,\\
Au&=&\lambda v,\\
 Lu&=&\lambda x,\\
 Ru&=&y.
\end{array}
\right.
\end{equation*}
Thus, $u=D_{\lambda^2}^{A,R}y$, and (3.4) follows. Further,
$D^{\mathbb{A} ,\mathbb{R} }_\lambda\in\mathcal{L}(\partial \mathbb{X}
,\mathbb{Z})$ for every Banach space
$\mathbb{Z}$ such that $D(\mathbb{A} ^\infty)\subset \mathbb{Z}\hookrightarrow
\mathbb{X} $. By 
 Lemma~3.5 the claim follows.
\end{proof}

\bigskip
We emphasize that the Dirichlet operators $D_\lambda^{\mathbb{A} ,\mathbb{R}
}$, $\lambda\in\rho (\mathbb{A}_0)$, are compact if and only if ${\partial X}$
is finite dimensional.

We now introduce a family of operators that 
will play an important role in the following.
By Lemma~3.8 and (3.1), we obtain the following.

\bigskip
{\bf Lemma~3.9.} {\sl Let $\lambda\in\rho(\mathbb{A}_0)$. Then the operator
$$\mathbb{B}_\lambda:=\tilde{\mathbb{B} }+\mathbb{B}  D_\lambda^{\mathbb{A} ,\mathbb{R} }\leqno(3.5)$$
exists, is represented by
$$\mathbb{B}_\lambda= B_1 D_{\lambda^2}^{A,R} +\left({1\over \lambda} B_3 +
B_4\right) L D_{\lambda^2}^{A,R},$$
and is bounded on ${\partial X}$.}

\bigskip
Using the family $(\mathbb{B}_\lambda)_{\lambda\in\rho(\mathbb{A}_0)}$ we can
now perform a useful factorization, similar to those discussed in 
\cite[\S~2]{[En99]}. This will allow us to investigate the spectral
properties of the matrix ${\mathcal{A}}$.

\bigskip
{\bf Lemma~3.10.} {\sl Let $\lambda\in\rho(\mathbb{A}_0)$.
Then the factorization
$$\lambda-{\mathcal{A}}=\mathcal{L}_\lambda {\mathcal{A}}_\lambda
\mathcal{M}_\lambda:= 
\begin{pmatrix}
I_{\mathbb{X} } & 0\\
-\mathbb{B} R(\lambda, \mathbb{A}_0) & I_{\partial \mathbb{X} }
\end{pmatrix}
\begin{pmatrix}
\lambda-\mathbb{A}_0 & 0\\
0 & \lambda-\mathbb{B}_\lambda
\end{pmatrix}
\begin{pmatrix}
I_{\mathbb{X} } & -D_\lambda^{\mathbb{A} ,\mathbb{R} }\\
0 & I_{\partial\mathbb{X} }
\end{pmatrix}\leqno(3.6)$$
holds, and for all $\mu\in{\bf C}$ we further have
$$\mu-{\mathcal{A}}=\mathcal{L}_\lambda
\begin{pmatrix}
\mu-\mathbb{A}_0 & 0\\
0 & \mu-\mathbb{B}_\lambda
\end{pmatrix}
\mathcal{M}_\lambda +(\mu-\lambda)
\begin{pmatrix}
0 & D_\lambda^{\mathbb{A} ,\mathbb{R} }\\ 
\mathbb{B}  R(\lambda,\mathbb{A}_0) &-\mathbb{B} R(\lambda,\mathbb{A}_0)D_\lambda^{\mathbb{A} ,\mathbb{R} }
\end{pmatrix}.\leqno(3.7)$$}

\begin{proof} Let $\lambda\in\rho(\mathbb{A}_0)$ and take
  ${\mathcal{U}}:={\mathbb{u} \choose \mathbb{v}}\in {\mathcal{X}}$. Observe
  first that ${\mathcal{U}}$ is in the domain 
  of the operator matrix $\mathcal{L}_\lambda {\mathcal{A}}_\lambda
  \mathcal{M}_\lambda$ if 
  and only if $\mathbb{u} -D_\lambda^{\mathbb{A} ,\mathbb{R}} \mathbb{v}\in
  D(\mathbb{A}_0)$, that is, if and only if 
$\mathbb{R} \left(\mathbb{u} -D_\lambda^{\mathbb{A} ,\mathbb{R} }\mathbb{v}\right)=\mathbb{R} \mathbb{u} -\mathbb{v}=0$. This shows that the
domains of the operators in (3.6) agree.
Moreover, we obtain
\begin{equation*}
\begin{array}{rl}
\begin{pmatrix}
  I_{\mathbb{X} } & 0\\
  -\mathbb{B} R(\lambda, \mathbb{A}_0) & I_{\partial \mathbb{X} }
\end{pmatrix}
&
\begin{pmatrix}
  \lambda-\mathbb{A}_0 & 0\\
  0 & \lambda-\mathbb{B}_\lambda
\end{pmatrix}
\begin{pmatrix}
  I_{\mathbb{X} } & -D_\lambda^{\mathbb{A} ,\mathbb{R} }\\
  0 & I_{\partial\mathbb{X}}
\end{pmatrix}
\begin{pmatrix}
  \mathbb{u} \\
  \mathbb{v}
\end{pmatrix}\\
=&
\begin{pmatrix}
  I_{\mathbb{X} } & 0\\
  -\mathbb{B} R(\lambda, \mathbb{A}_0) & I_{\partial\mathbb{X}}
\end{pmatrix}
\begin{pmatrix}
  \lambda-\mathbb{A}_0 & 0\\
  0 & \lambda-\mathbb{B}_\lambda
\end{pmatrix}
\begin{pmatrix}
  \mathbb{u} -D_\lambda^{\mathbb{A} ,\mathbb{R}}\mathbb{v}\\
  \mathbb{v}
\end{pmatrix}\\
=&
\begin{pmatrix}
  I_{\mathbb{X}} & 0\\
  -\mathbb{B} R(\lambda, \mathbb{A}_0) & I_{\partial\mathbb{X}}
\end{pmatrix}
\begin{pmatrix}
  (\lambda-\mathbb{A}_0) (\mathbb{u} -D_\lambda^{\mathbb{A},\mathbb{R}}
  \mathbb{v})\\
  \lambda \mathbb{v}-\mathbb{B}_\lambda\mathbb{v}
\end{pmatrix}\\
=&
\begin{pmatrix}
  (\lambda-\mathbb{A})(u-D_\lambda^{\mathbb{A},\mathbb{R}}\mathbb{v})\\
  -\mathbb{B} \mathbb{u}  +\mathbb{B} 
  D_\lambda^{\mathbb{A},\mathbb{R}}
  {\mathbb{v}+\lambda \mathbb{v}-\mathbb{B}_\lambda\mathbb{v}}
\end{pmatrix}
=\begin{pmatrix}
  \lambda-\mathbb{A}  & 0\\ 
  -\mathbb{B}  & \lambda-\tilde{\mathbb{B}}
\end{pmatrix}
\begin{pmatrix}
  \mathbb{u}\\
  \mathbb{v}
\end{pmatrix}
\end{array}
\end{equation*}
where we have used (3.5) and the fact that $D_\lambda^{\mathbb{A} ,\mathbb{R}
}$ maps $\partial\mathbb{X}  $ into $\ker(\lambda-\mathbb{A} )$, by
definition.  

To show (3.7), take $\mu\in{\bf C}$ and observe that 
\begin{equation*}
\begin{array}{rcl}
\mu-{\mathcal{A}} 
&=& 
(\mu-\lambda)I_{\mathcal{X}}+
\mathcal{L}_\lambda {\mathcal{A}}_\lambda \mathcal{M}_\lambda=
(\mu-\lambda)I_{\mathcal{X}}+\mathcal{L}_\lambda
\left[
\begin{pmatrix}
  \mu-\mathbb{A}_0 & 0\\
  0 & \mu-\mathbb{B}_\lambda
\end{pmatrix}
-(\mu-\lambda)I_{\mathcal{X}}
\right]
\mathcal{M}_\lambda\\
&=&\mathcal{L}_\lambda
\begin{pmatrix}
  \mu-\mathbb{A}_0 & 0\\
  0 & \mu-\mathbb{B}_\lambda
\end{pmatrix}
\mathcal{M}_\lambda+(\mu-\lambda)
(I_{\mathcal{X}}-\mathcal{L}_\lambda \mathcal{M}_\lambda).
\end{array}
\end{equation*}
One can check that
$$\mathcal{L}_\lambda\mathcal{M}_\lambda=
\begin{pmatrix}
I_{\mathbb{X} } & -D_\lambda^{\mathbb{A} ,\mathbb{R} }\\
-\mathbb{B}  R(\lambda,\mathbb{A}_0) & I_{\partial\mathbb{X}  }+
\mathbb{B} R(\lambda,\mathbb{A}_0)D_\lambda^{\mathbb{A} ,\mathbb{R} }
\end{pmatrix},\leqno(3.8)$$
and the claim follows.
\end{proof}

\bigskip
In many concrete cases, the spectrum of $A_0$, and hence by Lemma~3.4
of $\mathbb{A}_0$ are well-known. Hence it is interesting to decide whether a
given $\lambda\in\rho(\mathbb{A}_0)$ is a spectral value of the larger matrix
${\mathcal{A}}$. Using Lemma~3.10, we can now derive a partial characterization
whose main feature is the following: The spectrum and the
point spectrum (denoted by $\sigma$ and $P\sigma$, respectively) of a
$4\times 4$ operator matrix on $Y\times X\times{\partial X}\times{\partial X}$ is
characterized by means of the operator pencils
$(\mathbb{B}_\lambda)_{\lambda\in\rho(\mathbb{A}_0)}$ on ${\partial X}$.

\bigskip
{\bf Proposition~3.11.} {\sl For $\lambda\in\rho(\mathbb{A}_0)$ the equivalences
$$\lambda\in\sigma({\mathcal{A}})\iff\lambda\in\sigma(\mathbb{B}_\lambda)\qquad\hbox{and}
\qquad \lambda\in P\sigma({\mathcal{A}})\iff\lambda\in P\sigma(\mathbb{B}_\lambda)
\leqno(3.9)$$
hold. The set $\Gamma:=\left\{\lambda\in {\bf C}:
\lambda\in\rho(\mathbb{A}_0)\cap\rho(\mathbb{B}_\lambda)\right\}\subset\rho({\mathcal{A}})$ 
is nonempty, and for $\lambda\in\Gamma$ the resolvent operator of
${\mathcal{A}}$ is given by
$$R(\lambda,{\mathcal{A}})=
\begin{pmatrix}
R(\lambda,\mathbb{A}_0)+D_\lambda^{\mathbb{A} ,\mathbb{R} } R(\lambda,\mathbb{B}_\lambda)\mathbb{B} 
R(\lambda,\mathbb{A}_0) & D_\lambda^{\mathbb{A} ,\mathbb{R} } R(\lambda,\mathbb{B}_\lambda)\\
R(\lambda,\mathbb{B}_\lambda)\mathbb{B}  R(\lambda,\mathbb{A}_0) & R(\lambda,
\mathbb{B}_\lambda)
\end{pmatrix},\leqno(3.10)$$
where the entries are as in Lemmas~3.4, 3.7, and 3.8.}

\begin{proof}
Let $\lambda\in\rho(\mathbb{A}_0)$. Then the factorization (3.6)
holds. Observe that the operators $\mathcal{L}_\lambda$, $\mathcal{M}_\lambda$
are bounded and invertible, hence $\lambda-{\mathcal{A}}$ is invertible if and
only if the diagonal matrix ${\mathcal{A}}_\lambda$ is. We conclude that
$\lambda\in \sigma ({\mathcal{A}})$ if and only if
$\lambda\in\sigma(\mathbb{B}_\lambda)$. The latter equivalence in (3.9) follows
likewise. 

By Lemma~2.3(i) and Proposition~2.4, both the operators $\mathbb{A}_0$ and
${\mathcal{A}}$ are generators. Hence, their spectral bounds
$s(\mathbb{A}_0),s({\mathcal{A}})<\infty$. To show that $\Gamma\not=\emptyset$,
take thus $\lambda\geq\max\{s(\mathbb{A}_0),s({\mathcal{A}})\}$ and deduce by (3.9)
that $\lambda\in\Gamma$.

Finally, taking again into account (3.6) we obtain that for
$\lambda\in\Gamma$ there holds
$R(\lambda,{\mathcal{A}})=
\mathcal{M}_\lambda^{-1}{\mathcal{A}}_\lambda^{-1}\mathcal{L}_\lambda^{-1}$.
A direct computation now yields (3.10).
\end{proof}

\bigskip
Moreover, neglecting the trivial case of finite dimensional $X$, the
formula (3.10) allows us to obtain the following.

\bigskip
{\bf Theorem~3.12.} {\sl Let ${\rm dim}\;X=\infty$. Then the following
assertions are equivalent.

a) ${\mathcal{A}}$ has compact resolvent.

b) $\mathbb{A}_0$ has compact resolvent.

c) ${\partial X}$ is finite dimensional and the embeddings
$[D(A_0)]\hookrightarrow Y$ and $Y\hookrightarrow X$ are compact.}

\begin{proof} Take $\lambda\in\Gamma$ as defined in Proposition~3.11, so that
the resolvent operator $R(\lambda,{\mathcal{A}})$ is given by the formula
(3.10).

$a)\Rightarrow b)$ Let ${\mathcal{A}}$ have compact resolvent.
To begin with, the lower-right entry of
$R(\lambda,{\mathcal{A}})$ is the resolvent operator of a bounded operator on
${\partial X}$, thus it is compact if and only if ${\rm dim}\;
\partial\mathbb{X}  <\infty$. Therefore, the operator $D_\lambda^{\mathbb{A} ,\mathbb{R} }
R(\lambda,\mathbb{B}_\lambda)\mathbb{B}  R(\lambda,\mathbb{A}_0)$ is compact. Further, by
assumption the upper left entry of $R(\lambda,{\mathcal{A}})$, i.e.,
$R(\lambda,\mathbb{A}_0)+ D_\lambda^{\mathbb{A} ,\mathbb{R} } R(\lambda,\mathbb{B}_\lambda)\mathbb{B} 
R(\lambda,\mathbb{A}_0)$, is compact. It follows that their difference
$R(\lambda,\mathbb{A}_0)$ is compact.

$b)\Rightarrow a)$ Let now $\mathbb{A}_0$ have compact resolvent.
By Lemma~3.4, the identity on ${\partial X}$ is compact, and this implies that
${\rm dim}\;{\partial X}<\infty$. The claim now follows because the remaining
entries of $R(\lambda,{\mathcal{A}})$ are bounded operators with
finite-dimensional range.

$b)\Rightarrow c)$ In the following we consider $R(\lambda^2, A_0)$ as
a bounded, non-compact operator from $X$ to $[D(A_0)]$, and denote by
$i_Y$ and $i_X$ the embeddings $[D(A_0)]\hookrightarrow Y$ and
$Y\hookrightarrow X$, respectively. Let $\mathbb{A}_0$ have compact resolvent.
We have already seen that necessarily ${\rm
dim}\;{\partial X}<\infty$. Moreover, the $(1,1)$-entry of $R(\lambda,\mathbb{A}_0)$,
i.e., $i_Y \circ R(\lambda^2, A_0) \circ i_X$ is compact. Since also
the $(2,1)$-entry $A_0 R(\lambda^2, A_0) \circ i_X=\lambda^2 i_X \circ
i_Y R(\lambda^2, A_0) \circ i_X - i_X$ is compact, it follows that
$i_X$ is compact. Likewise, using the compactness of the (1,2)-entry of
$R(\lambda,\mathbb{A}_0)$, we can show the  compactness of $i_Y$, and the claim
is proven.

$c)\Rightarrow b)$ Finally, let the embeddings
$[D(A_0)]\hookrightarrow Y$ and $Y\hookrightarrow X$ be compact, and
${\rm dim}\;{\partial X}<\infty$. Then the embedding $[D(\mathbb{A}_0)]\hookrightarrow
\mathbb{X} $ is compact, and therefore $\mathbb{A}_0$ has compact resolvent.
\end{proof}

\bigskip
We now consider again the initial value problem associated with
(ABC) and prove a conjecture formulated in \cite[\S~5]{[GGG03]}. Prof. J. Goldstein has informed us that his student
C. Gal has recently obtained, by different methods, similar results on well-posedness and compactness issues. Gal's results  have been obtained simultaneously to, but independently of ours; they will appear in~\cite{[Ga04]}.

\bigskip
{\bf Corollary~3.13.} 
 {\sl Let the domain $\Omega$ be bounded. The matrix ${\mathcal{A}}$ associated with the
abstract version of $\rm(ABC)$ on $\Omega\subset{\mathbb R}^n$ has
compact resolvent if and only if $n=1$.}

\begin{proof} Recall that the embeddings 
$H^2(0,1)\hookrightarrow H^1(0,1)\hookrightarrow L^2(0,1)$
are compact. Then Theorem~3.12 yields the claim.
\end{proof}

\bigskip
To conclude this section, we mention that sharp results about the
essential spectrum $\sigma_{\rm ess}$ of the operator matrix arising
from the initial-boundary value problem associated with (ABC) have
been obtained in \cite[\S~3]{[Be76]}. The proofs therein are very technical,
and only work if the domain $\Omega$ is bounded.

The formula (3.7) can however be used to obtain some results
about $\sigma_{\rm ess}({\mathcal{A}})$, too. Our Propositions~3.14 and 4.1
below complement the results due to Beale. In particular,
Proposition~3.14 also applies if we consider the motivating equation 
(ABC) to take place on the unbounded domain $\Omega={\bf R}_+$.

\bigskip
{\bf Proposition~3.14.} {\sl Let ${\partial X}$ be finite dimensional. Then the
essential spectrum of ${\mathcal{A}}$ is given by 
$$\sigma_{\rm ess}({\mathcal{A}})= \sigma_{\rm ess}(\mathbb{A}_0),$$
and for the Fredholm index we have
$$ind({\mathcal{A}}-\mu)=ind(\mathbb{A}_0-\mu)\qquad\hbox{for all}\;
\mu\not\in \sigma_{\rm ess}(\mathbb{A}_0).$$}

\begin{proof} To begin with, we recall that the essential spectrum does
neither change under compact additive perturbations, nor under
similarity transformations (i.e., bounded invertible multiplicative perturbations).

Fix $\lambda\in\rho(\mathbb{A}_0)$, take into account (3.8), and observe
that $I_{\mathcal{X}}-\mathcal{L}_\lambda \mathcal{M}_\lambda$ is a compact
operator on ${\mathcal{X}}$. Moreover, $\mathcal{L}_\lambda$,
$\mathcal{M}_\lambda$ are bounded and invertible. Thus, to decide whether a
given $\mu\in{\bf C}$ is in the essential spectrum of ${\mathcal{A}}$, by
(3.7) it suffices to check whether 0 is in the essential spectrum of the
operator matrix 
$$\begin{pmatrix}
\mu-\mathbb{A}_0 & 0\\
0 & \mu-\mathbb{B}_\lambda
\end{pmatrix}=
\begin{pmatrix}
\mu-\mathbb{A}_0 & 0\\
0 & 0
\end{pmatrix}
+
\begin{pmatrix}
0 & 0\\ 
0 & \mu-\mathbb{B}_\lambda
\end{pmatrix}.$$
The second addend is a bounded operator with finite dimensional range,
hence it does not affect the essential spectrum of the operator matrix
on the left-hand side, and the claim follows.
\end{proof}

\section{The special case of {$B_3=0$}}

After setting $y=\dot{x}$, $\rm(AIBVP_2)$ can
equivalently be written as the second order problem with {\sl
integro-differential boundary conditions}
\begin{equation*}
\left\{
\begin{array}{rcll}
 \ddot{u}(t)&=& Au(t), &t\in {\bf R},\\
 \dot{y}(t)&=& B_1u(t)+B_2\dot{u}(t)+ B_3 \left(h+\!\int_0^t\!
 y(s)\;ds\right) + B_4 y(t), &t\in {\bf R},\\
  y(t)&=& Lu(t), &t\in {\bf R},\\
  u(0)&=&f, \qquad \dot{u}(0)=g,&\\
  y(0)&=&j.&
\end{array}
\right.
\end{equation*} 
In the special case of $B_3=0$, which we assume throughout this
section, the initial value $x(0)=h$ is therefore superfluous, and we
obtain an abstract second  order problem with first order dynamical
boundary conditions. Similar problems have been discussed, among
others, in \cite{CENP05}, and in fact some well-posedness result therein,
cf. \cite[Thm.~2.2]{CENP05}, can be interpreted as a corollary of our 
Theorem~2.7. Moreover, observe that we can now replace $\mathbb{X} =Y\times
X\times {\partial X}$ by $\mathbb{X} =Y\times X$, and the operator matrix $\mathbb{A} $ as
defined in (3.1) by 
$$\mathbb{A} =
\begin{pmatrix}
0 & I\\ 
A & 0
\end{pmatrix},\qquad
D(\mathbb{A} )=D(A)\times Y.\leqno(4.1)$$ 
Accordingly, the operators $\mathbb{R} $ and $\mathbb{B} $ become
$$\mathbb{R} =
\begin{pmatrix}
R & 0
\end{pmatrix},\qquad D(\mathbb{R} )= D(A)\times X,\qquad\qquad
\mathbb{B} =
\begin{pmatrix}
B_1+B_4 B_2 & 0
\end{pmatrix},\qquad D(\mathbb{B} )=\mathbb{X} .$$
Then the operator matrix ${\mathcal{A}}$ defined in (2.8) becomes
$${\mathcal{A}}=
\begin{pmatrix}
0 & I & 0\\ 
A & 0 & 0\\ 
B_1+B_4 B_2 & 0 & B_4
\end{pmatrix},\qquad
D({\mathcal{A}})=\left\{
\begin{pmatrix}
u\\ v\\ y
\end{pmatrix}
\in D(A)\times Y\times {\partial X}\:
Lu=B_2u\right\}.$$
The main difference with the general setting of Section~3 is that
the resolvent of $\mathbb{A}_0$ as well as the Dirichlet operators associated
with $\mathbb{A} $ and $\mathbb{R} $ can be compact also in the case of ${\rm
dim}\;\partial X=\infty$.

This allows us to obtain the following.

\bigskip
{\bf Proposition~4.1.} {\sl Assume that there exists a Banach space
$Z$ containing $D(A)$ and compactly embedded in $Y$, and that also
the embedding of $Y$ in $X$ is compact. Then  
$$\sigma_{\rm ess}({\mathcal{A}})=\sigma_{\rm ess}(B_4)\qquad\hbox{and}\qquad
ind({\mathcal{A}}-\lambda)= ind(B_4-\lambda)\quad\hbox{for}\;\lambda\not\in
\sigma_{\rm ess}(B_4).$$
In particular, $\sigma_{\rm ess}({\mathcal{A}})=\emptyset$
if and only if ${\partial X}$ is finite dimensional.}

\begin{proof} Fix $\lambda\in\rho(\mathbb{A}_0)$. Then for all $\mu\in{\bf C}$ the
factorization (3.7) holds. Taking into account Remark~3.3(b), we see
that by assumption $D_\lambda^{\mathbb{A} ,\mathbb{R} }$ and $\mathbb{B}  R(\lambda, \mathbb{A}_0)$ are
compact operators from $\mathbb{X} $ to $\partial\mathbb{X}  $ and from $\partial\mathbb{X}$ to
$\mathbb{X}$, respectively. Thus, reasoning as in the proof of
Proposition~3.14, we obtain that 
$\mu\in\sigma_{\rm ess}({\mathcal{A}})$ if and only if $\mu\in\sigma_{\rm
ess}(B_4)$. Here we have used the fact that $\sigma_{\rm
ess}\big(\mathbb{B}_\lambda)=\sigma_{\rm ess}(\tilde{\mathbb{B} })=\sigma_{\rm
ess}(B_4)$, and that $D(\mathbb{A}_0)$ is
compactly embedded in $\mathbb{X} $, i.e., $\mathbb{A}_0$ has empty essential spectrum.

Finally, recall that a bounded operator has empty essential spectrum
if and only if it acts on a finite dimensional space,
cf. \cite[\S~IV.1.20]{[EN00]}.
\end{proof}

\bigskip 
{\bf Example~4.2.} In the context of our motivating equation (ABC),
the assumption $B_3=0$ means that $k\equiv 0$, hence the
initial-boundary value problem becomes 
\begin{equation}\tag{*}
\left\{
\begin{array}{rcll}
 \ddot{\phi}(t,x)&=& c^2 \Delta \phi(t,x), &t\in {\bf R},\;
 x\in \Omega,\\
 \ddot{\delta}(t,z)&=&  -{d(z)\over m(z)}\dot{\delta}(t,z)
-{\rho(z)\over m(z)}\dot{\phi}(t,z), &t\in {\bf  R},\;z\in\partial \Omega,\\
  \dot{\delta}(t,z)&=& {\partial\phi \over\partial\nu}(t,z),
  &t\in {\bf R},\;z\in\partial \Omega,\\
\phi(0,\cdot)&=&f, \qquad\dot{\phi}(0,\cdot)=g,&\\
\dot{\delta}(0,\cdot)&=&j,
\end{array}
\right.   
\end{equation}
on a bounded open domain $\Omega\subset{\bf R}^n$. 
Observe that $D(A)\subset Z:= H^{3\over 2}(\Omega)$, and for
$Y=H^1(\Omega)$ the embeddings $H^{3\over 2}(\Omega)\hookrightarrow
H^1(\Omega) \hookrightarrow L^2(\Omega)$ are compact by
\cite[Vol.~I, Thm.~1.16.1]{[LM72]}.
Thus, by Proposition~4.1 the essential spectrum of
${\mathcal{A}}$ agrees with the essential spectrum of the bounded
multiplication operator
$$(B_4 u)(z)= -{d(z)\over m(z)}u(z),\qquad u\in L^2(\partial \Omega),\; z\in\partial \Omega.$$
The essential spectrum of $B_4$ cannot be empty unless ${\partial X}$ is finite
dimensional, thus the essential spectrum of ${\mathcal{A}}$ cannot be
empty unless $n=1$.\qed

\bigskip
Recall that in the context of our motivating equation (ABC) we always
have $B_1=0$. However, if $B_3=0$ and the feedback $B_1$ is instead of
the form $B_1=-B_4 B_2$, then we obtain sharper spectral results. In
fact, in this case ${\mathcal{A}}$ becomes a diagonal block matrix (with
nondiagonal domain). Observe that the resolvent 
  set of $\mathbb{A}_0$ is $\left\{\lambda\in {\bf C}\:
\lambda^2\in\rho(A_0)\right\}$, instead of $\left\{0\not=\lambda\in
{\bf C}\: \lambda^2\in\rho(A_0)\right\}$ as in Lemma~3.4, which is
remarkable because in the context of our motivating equation (ABC)
$A_0$ is the Laplacian with Robin boundary conditions, so that
$0\in\rho(A_0)$ and therefore also $0\in\rho(\mathbb{A}_0)$.

\bigskip
We can now obtain an easier version of the equivalence (3.9) --
getting rid of the operator pencil $\mathbb{B}_\lambda$, $\lambda\in\rho(\mathbb{A}_0)$
 -- and derive an alternative characterization of the spectrum of ${\mathcal{A}}$
 that compliments the one already obtained in Proposition~3.11.

\bigskip
{\bf Corollary~4.3.} {\sl Let $B_3=0$ and $B_1=-B_4 B_2$. Then the
following hold.

i) If $\lambda^2\in\rho(A_0)$, then
$$\lambda\in\sigma({\mathcal{A}})\iff
\lambda\in\sigma(B_4)\qquad\hbox{and}
\qquad \lambda\in P\sigma({\mathcal{A}})\iff\lambda\in P\sigma(B_4).$$
\qquad If $\lambda^2\in\rho(A_0)$ and $\lambda\in\rho(B_4)$, then the
resolvent operator $R(\lambda,{\mathcal{A}})$ is given by
$$R(\lambda,{\mathcal{A}}):=
\begin{pmatrix}
\lambda R(\lambda^2,A_0) & R(\lambda^2,A_0) & D_{\lambda^2}^{A,R}
R(\lambda,B_4)\\ 
A_0 R(\lambda^2,A_0) & \lambda R(\lambda^2,A_0) & \lambda D_{\lambda^2}^{A,R} R(\lambda,B_4)\\
0 & 0 & R(\lambda,B_4)
\end{pmatrix}.$$

ii) If $\lambda\not\in P\sigma(B_4)$, then
$$\lambda\in P\sigma({\mathcal{A}})\iff \lambda^2\in P\sigma(A_0).$$}

\begin{proof} i) Lemma~3.9 yields that
$$\mathbb{B}_\lambda= (B_1+B_4L)D_{\lambda^2}^{A,R},\qquad
\lambda\in\rho(\mathbb{A}_0).\leqno(4.2)$$
Now, taking into account (3.1) and Proposition~3.11 the claim follows.

ii) Let
$$({\mathcal{A}}-\lambda){\mathcal{U}}=
\begin{pmatrix}
v-\lambda u\\
Au-\lambda v\\
(B_4-\lambda)Ru
\end{pmatrix}=0.$$
Thus, we obtain that $(A_0-\lambda^2)u=0$ and the claim follows.
\end{proof}

\bigskip
Hence, we can sometimes obtain a complete characterization of
the point spectrum.

\bigskip
{\bf Corollary~4.4.} {\sl Let $B_3=0$ and $B_1=-B_4 B_2$. Assume that
$$\{\lambda\in P\sigma(B_4)\:
\lambda^2\in \sigma(A_0)\}=\emptyset.\leqno(4.4)$$
Then
$$P\sigma({\mathcal{A}})=\{\lambda\in {\bf C}\: \lambda^2\in P\sigma(A_0)
\;\;\hbox{or}\;\; \lambda\in P\sigma(B_4)\}.$$}

\medskip
Observe that the condition (4.4) is in particular satisfied whenever
$A_0$ is self-adjoint and invertible and $B_4$ has no eigenvalues on
$i{\bf R}\setminus\{0\}$.

\bigskip 
{\bf Example~4.5.} Let $k\equiv0$, and $B_1=-B_4 B_2$, that is,
$$(B_1 f)(z)=-{\rho(z) d(z)\over m^2(z)}f(z),\qquad f\in
H^1(\Omega),\; z\in\partial \Omega.$$
Hence, we revisit Example~4.2 and consider a version of $(*)$ where 
we replace the second equation by
$$\ddot{\delta}(t,z)=
 -{d(z)\over m(z)}\dot{\delta}(t,z)-{\rho(z) d(z)\over m^2(z)}\phi(z)
-{\rho(z)\over m(z)}  \dot{\phi}(t,z),
\qquad t\in {\bf  R},\;z\in\partial \Omega.$$
The Laplacian with Robin boundary conditions is self-adjoint and
injective (with compact resolvent), thus the spectrum of $A_0$
consists of  countably many strictly negative values diverging to $-\infty$,
cf. \cite[\S~IV.3--4]{[Mi78]}. On the other hand, $B_4$ is a multiplication
operator. Thus, its spectrum agrees with the essential range of
the function $-{d(\cdot)\over m(\cdot)}$, while its point spectrum 
is given by $\big\{\lambda\in {\bf C}\: \mu\{z\in \partial \Omega\:
d(z)+\lambda m(z)=0\}>0\big\}$, cf. \cite[Ex.~I.4.13(8)]{[EN00]}.

In particular, since by assumption $m$ is
real valued, the condition (4.4) is satisfied if the essential range
of $d$ does not contain any point on $i{\bf R}$. In this case by
Corollary~4.3 we can completely characterize the point 
spectrum of the matrix ${\mathcal{A}}$ associated with $(*)$. In addition,
${\mathcal{A}}$ turns out to be invertible on ${\mathcal{X}}$. (Taking into
account (4.2), one can see that this is not the case if $B_1=0$,
because then $0\in\rho({\mathcal{A}})$ if and only if $0\in\rho(B_4 L
D^{A,R}_0)$. 
This does not hold, since the normal derivative $L$ vanishes on
the set of constants.)

\section{Neutral acoustic boundary conditions}

Among the so-called {\sl boundary contact problems} discussed by
Belinsky in \cite[\S~3]{[Bl00]}, the Timoschenko model
\begin{equation}\tag{TM}
\left\{
\begin{array}{rcll}
 \ddot{\phi}(t,x)&=& c^2 \Delta \phi(t,x,z), &t\in {\bf R},\;
 x\in \Omega,\\
 {\partial\phi \over\partial\nu}(t,z)&=&0, &t\in {\bf R},
    \;z\in\Gamma_0,\\
 m(1-\Delta) \ddot{\delta}(t,z)&=&
 -d(z)\dot{\delta}(t,z)-k(z){\delta}(t,z)-\rho
\dot{\phi}(t,z), &t\in {\bf R},\;z\in\Gamma_1,\\
  \dot{\delta}(t,z)&=&{\partial\phi \over\partial\nu}(t,z),
  &t\in {\bf R},\;z\in\Gamma_1,
\end{array}
\right.
\end{equation}
is particularly interesting, because it can be seen as a wave 
 equation equipped with {\sl neutral} acoustic boundary
 conditions. The aim of this section is to show how the methods introduced
 above can be applied to the present situation with minor changes.

For the geophysical explanation of this model we refer the reader to
 \cite{[Bl00]}. We only mention that the system (TM) models an
 ocean waveguide  $\Omega$ covered (on the part $\Gamma_1$ of his
 surface $\partial \Omega$) by a thin pack ice layer with inertia of
 rotation. Belinsky investigates such a system for $\Omega\subset {\bf
 R}^2$ only and obtains some spectral properties. 

Here the boundary $\partial \Omega$ is the disjoint union of
 $\Gamma_0,\Gamma_1$. Observe that, due to technical reasons, we
 consider the case of a medium of {\sl homogeneous} 
 density $\rho$ filling a domain whose boundary has {\sl homogeneous}
 mass $m$. However, we still allow $k$ and $d$ to be essentially
 bounded functions, whereas Belinsky assumes them to be constant.

\bigskip
In this section we cast such an equation into an abstract framework, and
discuss its well-posedness -- this is a new result to our knowledge.

To begin with, we introduce an operator $M$ that will appear in the new
neutral acoustic boundary conditions.

\bigskip 
\goodbreak
{\bf Assumption~5.1}

(A$_8$) The operator $M:D(M)\subset {\partial X}\to {\partial X}$ is linear,
closed, and satisfies $1\in\rho(M)$.

\bigskip
We can now consider the abstract second order initial-boundary value
problem obtained by replacing the second equation in $\rm(AIBVP_2)$ by
$$\ddot{x}(t)-M\ddot{x}(t)= B_1u(t)+B_2\dot{u}(t)+ B_3 x(t)
+ B_4\dot{x}(t),\qquad t\in {\bf R}.$$
Thus, our aim is to show the well-posedness of the problem
\begin{equation}\tag{AIBVP$^\diamond_2$}
\left\{
\begin{array}{rcll}
 \ddot{u}(t)&=& Au(t), &t\in {\bf R},\\
 \ddot{x}(t)&=& B^\diamond_1u(t)+B^\diamond_2\dot{u}(t)+ B^\diamond_3 x(t) +
 B^\diamond_4\dot{x}(t),
 &t\in {\bf R},\\
  \dot{x}(t)&=& Lu(t), &t\in {\bf R},\\
  u(0)&=&f,\quad \dot{u}(0)=g,&\\
  x(0)&=&h,\quad \dot{x}(0)=j,&
\end{array}
\right.
\end{equation}
on $X$ and ${\partial X}$, where
$${B^\diamond_i}:= R(1,M)B_i, \qquad i=1,2,\leqno(5.1)$$
are bounded operators from $Y$ to ${\partial X}$, and
$${B^\diamond_i}:= R(1,M)B_i, \qquad i=3,4,\leqno(5.2)$$
are bounded operators on ${\partial X}$. Similarly, we consider the operator 
$${R^\diamond}:=L-B^\diamond_2.\leqno(5.3)$$
Observe now that, after replacing $R$ by $R^\diamond$ and $B_i$ by $B^\diamond_i$,
$i=1,2,3,4$, all the corresponding Assumptions~2.1 are satisfied,
except for (A$_3$) and (A$_7$). To fill this gap, we replace them by
the following.

\bigskip 
{\bf Assumptions~5.2} 

(A$'_3$) The operator $R^\diamond:D(A)\to {\partial X}$ is linear and
surjective.

(A$'_7$) The restriction $A^\diamond_0:=A_{\arrowvert\!\ker R^\diamond}$
generates 
a cosine operator function with associated phase space $Y\times X$.

\bigskip
Under the assumptions (A$_1$), (A$_2$), (A$'_3$), (A$_4$)--(A$_6$),
(A$'_7$), and (A$_8$) we promptly obtain the main result of this
section. 

\bigskip
{\bf Proposition~5.3.} {\sl The problem (AIBVP$^\diamond_2$) with abstract
neutral acoustic boundary conditions is well-posed.}

\begin{proof} Consider the operator matrix $\mathbb{A} $ as introduced in
  (2.2) and define the operators
$$\mathbb{R}^\diamond:=
\begin{pmatrix}
  R^\diamond & 0 & 0
\end{pmatrix},
\qquad
D(\mathbb{R}^\diamond):=D(\mathbb{A} ),$$
$$\mathbb{B}^\diamond:=
\begin{pmatrix}
B^\diamond_1+B^\diamond_4 B^\diamond_2 & 0 & B^\diamond_3
\end{pmatrix},
\qquad 
D(\mathbb{B}^\diamond):=\mathbb{X},$$
$$\tilde{\mathbb{B}}^\diamond:=B^\diamond_4,
\qquad  D(\tilde{\mathbb{B} }^\diamond):=\partial \mathbb{X} .$$
We can now directly check that properties analogous to those in
Lemma~2.3 are satisfied. Therefore, the well-posedness of
\begin{equation}\tag{$\mathbb{AIBVP}^\diamond$}
\left\{
\begin{array}{rcll}
 \dot{\mathbb{u} }(t)&=& \mathbb{A} \mathbb{u} (t), &t\in {\bf R},\\
 \dot{\mathbb{x} }(t)&=& \mathbb{B}^\diamond\mathbb{u} (t)+\tilde{\mathbb{B}
 }^\diamond\mathbb{x} (t), &t\in {\bf  R},\\ 
 \mathbb{x} (t)&=&\mathbb{R}^\diamond\mathbb{u} (t),  &t\in {\bf R},\\
  \mathbb{u} (0)&=&\mathbb{u}_0,&\\
  \mathbb{x} (0)&=&\mathbb{x}_0,&
\end{array}
\right.
\end{equation}
follows like in Proposition~2.4.
Finally, reasoning like in Lemma~2.6 one obtains the equivalence
 between $(\mathbb{AIBVP}^\diamond)$ and $\rm(AIBVP_2^\diamond)$, and the claim
 follows.
\end{proof}

\bigskip
It is usually not easy to check directly whether the assumption (A$'_3$) is
satisfied. Thus, we propose some conditions on interpolation
inequalities that ensure its validity: In many concrete cases the restriction 
$$A_L:=A_{\arrowvert\!\ker L}$$ 
generates an analytic semigroup on $X$, and we can therefore consider
the interpolation spaces 
$$X_\theta:=[D(A_L),X]_{(1-\theta)},\qquad 0\le \theta\le 1,$$
cf. \cite[Chapt.~1]{[Lu95]} for the abstract theory and \cite[Vol.~I,  Chapt.~1]{[LM72]} for concrete spaces.

\bigskip
{\bf Lemma~5.4.} {\sl Let $L$ be surjective and $A_L$ generate an
analytic semigroup. Assume that for some $0<\alpha<1$ and
$\tilde\epsilon\in (0,1-\alpha)$ one has $D(A)\subset
X_{\alpha+\tilde\epsilon}$ and $X_{\alpha}\hookrightarrow Y$. Then
$R^\diamond$ is surjective.}

\begin{proof} First of all we remark the existence of bounded Dirichlet
operators $D_\lambda^{A,L}$ associated with $A$ and $L$,
$\lambda\in\rho(A_L)$. Further, we can reason as in \cite[Lemma~2.4]{[GK91]} 
and obtain that 
$$\| D_\lambda^{A,L}\|_{_{{\partial X},Y}}=
O(\mid\!\lambda\!\mid^{-\tilde\epsilon})\qquad\hbox{for}\;\;
\mid\!\lambda\!\mid\to\infty,\; {\rm Re}(\lambda)>0.$$
Thus, we can choose $\lambda_0$ large enough so that
$\| B^\diamond_2 D_{\lambda_0}^{A,L}\|< 1$, 
and consequently we can invert the operator $I_{\partial X}- B^\diamond_2
D_{\lambda_0}^{A,L}=R^\diamond D_{\lambda_0}^{A,L}$.
y
We are now in the position to prove the surjectivity of $R^\diamond$.
Take $x\in{\partial X}$, and observe that for the vector
$$u:=D_{\lambda_0}^{A,L}\left(R^\diamond D_{\lambda_0}^{A,L}\right)^{-1}x$$
there holds $R^\diamond u=x$.
\end{proof}

\bigskip
We revisit the motivating example.

\bigskip
{\bf Theorem~5.5.} {\sl The initial value problem associated
with the wave equation with neutral acoustic boundary conditions
$\rm(TM)$ on an open bounded domain $\Omega\subset {\bf R}^n$,
$n\geq1$, with smooth boundary $\partial \Omega=\Gamma_0\cup \Gamma_1$ is
well-posed. In particular, for all initial data 
$$\phi(0,\cdot)\in H^2(\Omega),\qquad
\dot{\phi}(0,\cdot)\in H^1(\Omega),\qquad \delta(0,\cdot)\in
L^2(\Gamma_1),$$
$$\hbox{and}\qquad\dot{\delta}(0,\cdot)\in L^2(\Gamma_1)\qquad\hbox{such
that}\qquad{\partial \phi\over \partial
\nu}(0,\cdot)=\dot{\delta}(0,\cdot)$$
there exists a classical solution on
$(H^1(\Omega),L^2(\Omega),L^2(\Gamma_1))$ continuously depending on
them.}

\begin{proof} We take over the setting introduced in Theorem~2.7 and adapt it
to the current problem. We let 
$$X:=L^2(\Omega),\qquad
Y:=H^1(\Omega),\qquad {\partial X}:=L^2(\Gamma_1).$$
Moreover, we set
$$A:=c^2\Delta,\qquad D(A):=\left\{u\in H^{3\over
2}(\Omega)\: \Delta u\in L^2(\Omega),\; {\partial
u\over\partial\nu}_{\big\arrowvert\! \Gamma_0}\!\!=0\right\},$$
$$(Rf)(z)={\partial f\over\partial \nu}(z)+{\rho(z)\over
m(z)} f(z),\qquad f\in D(R)=D(A),\; z\in\Gamma_1,$$
$$B_1=0,\qquad (B_2 f)(z):=-{\rho(z)\over m(z)} f(z),\qquad f\in
H^1(\Omega),\; z\in\Gamma_1,$$
$$(B_3 g)(z):=-{k(z)\over m(z)}g(z),\qquad (B_4
g)(z):=-{d(z)\over m(z)}g(z),\qquad g\in L^2(\Gamma_1),\;
z\in\Gamma_1.$$
Further, we introduce the operator 
$$M:=\Delta,\qquad D(M):=H^2(\Gamma_1),$$
that is, the Laplace--Beltrami operator on $\Gamma_1$, and we can now
define the auxiliary operators $B^\diamond_i$, $i=1,2,3,4$, and $R^\diamond$
like in (5.1)--(5.3).

The operator $M$ clearly satisfies the assumption (A$_8$), 
hence only the assumptions (A$'_3$) and (A$'_7$) still need to
be checked.

To check (A$'_3$), first observe that $L=R+B_2$ is the normal
derivative on $\Gamma_1$, which is surjective by \cite[Vol.~I,
  Thm.~2.7.4]{[LM72]}. Moreover, consider the operator 
$A_L=A_{\arrowvert\!\ker L}$ and observe that 
$$D(A_L)=\left\{u\in H^2(\Omega)\:
{\partial u\over\partial
\nu}_{\big\arrowvert\!\partial \Omega}\!\!=0\right\},$$ 
that is, $A_L$ turns out to be the Neumann Laplacian. Thus, $A_L$
generates an analytic
semigroup on $X$. Further, by \cite[Vol.~II, (4.14.32)]{[LM72]}, we obtain that
$$\left[ [D(A_L)], L^2(\Omega)\right]_\theta =
H^{2(1-\theta)}(\Omega)\qquad\hbox{for
all}\; \theta\in\left({1\over 4},1\right].$$ 
Therefore, for $X_\alpha=\left[ [D(A_L)],
L^2(\Omega)\right]_{1-\alpha}$, we obtain that
$$D(A)\subset H^{3\over 2}(\Omega)\subset X_{\alpha+\epsilon}
\subset X_{\alpha}=H^1(\Omega)=Y\qquad\hbox{for}\;
\alpha={1\over 2}\;\hbox{and all}\; 0\le\epsilon<{1\over 4}.$$
By Lemma~5.4 the assumption~(A$'_3$) is thus checked.

Finally, to check (A$'_7$) we show that $A^\diamond_0$ is self-adjoint
and dissipative. Recall that we are assuming $\rho$ and $m$ to be
real constants, and observe that  
$$D(A^\diamond_0)=\left\{f\in H^2(\Omega)\: 
{\partial u\over\partial \nu}_{\big\arrowvert\Gamma_0}\!\!=0,\;
\left({\partial u\over\partial
\nu}_{\big\arrowvert\Gamma_1}\!\!+{\rho\over 
m}R(1,M)(u_{\arrowvert\Gamma_1})\!\!\right)=0\right\}.$$ 
Take $u,v\in D(A^\diamond_0)$ and obtain that
\begin{equation*}
\begin{array}{rcl}
\big<A^\diamond_0 u,v\big>_X
&=&
\int_\Omega \Delta u\cdot\overline{v}\;dx= 
\int_{\partial \Omega} {\partial u\over
\partial\nu}\cdot\overline{v}\;d\sigma - 
\int_\Omega \nabla u\cdot\overline{\nabla v}\;dx\\
&=&-{\rho\over m}\int_{\Gamma_1}
R(1,M)(u_{\arrowvert\Gamma_1})\cdot\overline{v}\;d\sigma - 
\int_\Omega \nabla u\cdot\overline{\nabla v}\;dx.
\end{array}
\end{equation*}
Taking into account the positivity and the self-adjointness of
the operator $R(1,M)$, one obtains that $A^\diamond_0$ is self-adjoint and
dissipative, hence it generates a cosine operator function. We still
need to check that its phase space actually is $Y\times X$. To 
do so, it is convenient to use a variational argument.
Integrating by parts one sees that
$A^\diamond_0$ is {\sl not} invertible, hence we need to consider its
(invertible) perturbation $\tilde{A}:=(A^\diamond_0-I_X)$, which also
generates a cosine operator function. Observe that, due to the
boundedness of the perturbation, $\tilde{A}$ and $A^\diamond_0$ have same form
domain as well as same phase space. Reasoning as in the last lines
of Example~2.8, the claim therefore follows if we can show that the
form domain of $\tilde{A}$ is $Y$. 

In fact, the sesquilinear form associated with $\tilde{A}$ is
$$a(u,v):= \int_\Omega \nabla u\cdot\overline{\nabla v}\;dx 
+ {\rho\over m}\int_{\Gamma_1} 
R(1,M)^{1\over 2}(u_{\arrowvert\Gamma_1})\cdot
\overline{R(1,M)^{1\over 2}(v_{\arrowvert\Gamma_1}})\;d\sigma
-\int_\Omega u\cdot\overline{v}\;dx,$$
whose domain is actually $H^1(\Omega)=Y$.
\end{proof}
 
\bigskip
\noindent
{\bf Acknowledgment.} I gratefully thank Prof. Jerry Goldstein
(University of Memphis) for bringing the theory of acoustic boundary
conditions to my attention.


\begin{thebibliography}{10}

\bibitem{[ABHN01]} W. Arendt, C.J.K. Batty, M. Hieber, and F.
Neubrander, \emph{Vector-valued Laplace Transforms and Cauchy
Problems}, Monographs in Mathematics \textbf{96}, Birkh\"auser Verlag 2001.

\bibitem{[Be76]} J.T. Beale, \emph{Spectral properties of an acoustic
boundary condition}, Indiana Univ. Math. J. \textbf{25} (1976), 895--917.

\bibitem{[BR74]} J.T. Beale and S.I. Rosencrans, \emph{Acoustic
boundary conditions}, Bull. Amer. Math. Soc. \textbf{80} (1974),
1276--1278.

\bibitem{[Bl00]} B.P. Belinsky, \emph{Wave propagation in the ice-covered
ocean wave guide and operator polynomials}, in: H.G.W. Begehr,
R.P. Gilbert, and J. Kajiwara (eds.), ``Proceedings of the Second
ISAAC Congress'', Kluwer Academic Publishers 2000, 1319--1333.

\bibitem{[CENN03]} V. Casarino, K.-J. Engel, R. Nagel, and G. Nickel,
\emph{A semigroup approach to boundary feedback systems}, Integral
Equations Oper. Theory {\bf 47} (2003), 289--306.

\bibitem{CENP05} V. Casarino, K.-J. Engel, G. Nickel, and S. Piazzera, \emph{Decoupling techniques for wave equations with dynamic boundary conditions}, Disc. Cont. Dyn. Syst {\bf 12} (2005), 761-772.

\bibitem{[DL90]} R. Dautray and J.-L. Lions, \emph{Mathematical
Analysis and Numerical Methods for Science and Technology. Voll.
1--2}, Springer-Verlag 1988-1990.

\bibitem{[En99]} K.-J. Engel, \emph{Spectral theory and generator
property for one-sided coupled operator matrices}, Semigroup Forum
\textbf{58} (1999), 267--295.

\bibitem{[EN00]} K.-J. Engel and R. Nagel, \emph{One-Parameter Semigroups
for Linear Evolution Equations}, Graduate Texts in Mathematics \textbf{194},
Springer-Verlag 2000. 

\bibitem{[Fa85]} H.O. Fattorini, \emph{Second Order Linear Differential
Equations in Banach Spaces}, Mathematics Studies \textbf{108},
North-Holland 1985.

\bibitem{[Ga04]} C. Gal, Ph.D. thesis (in preparation).

\bibitem{[GGG03]} C. Gal, G.R. Goldstein, and J.A. Goldstein, \emph{Oscillatory boundary conditions for acoustic wave equations}, J. Evol. Equations {\bf 3} (2004), 623--636.

\bibitem {[GK91]} G. Greiner and K. Kuhn,
\emph{Linear and semilinear boundary conditions: the analytic case},
in: Ph. Cl\'ement, E. Mitidieri, and B. de Pagter (eds.), ``Semigroup
Theory and Evolution Equations'' (Proceedings Delft 1989),  Lecture
Notes in Pure and Appl. Math. \textbf{135}, Marcel Dekker 1991, 193--211.

\bibitem{[Ka95]} T. Kato, \emph{Perturbation Theory for Linear Operators},
Classics in Mathematics, Springer-Verlag 1995.

\bibitem{[KMN03]} M. Kramar, D. Mugnolo, and R. Nagel, \emph{Semigroups
 for initial-boundary value problems}, in: M. Iannelli and G. Lumer
 (eds.): ``Evolution Equations 2000: Applications to Physics,
 Industry, Life Sciences and Economics'' (Proceedings Levico Terme
 2000), Progress in Nonlinear Differential Equations, Birkh\"auser
 2003, 277--297.

\bibitem{[Kr61]} V.N. Krasil'nikov, \emph{On the solution of some
 boundary-contact problems of linear hydrodynamics},
J. Appl. Math. Mech. \textbf{25} (1961), 1134--1141.

\bibitem{[LM72]} J.L. Lions and E. Magenes, \emph{Non-Homogeneous
Boundary Value Problems and Applications. voll. I--II}, Grundlehren der
mathematischen Wissenschaften \textbf{181--182}, Springer-Verlag 1972.

\bibitem{[Lu95]} A. Lunardi, \emph{Analytic Semigroups and Optimal 
Regularity in Parabolic Problems}, Progress in Nonlinear Differential 
Equations and their Applications \textbf{16}, Birkh\"auser 1995.

\bibitem{[Mi78]} V.P. Mikhajlov, \emph{Partial Differential Equations},
Mir Publishers 1978.

\bibitem{[MI68]} P.M. Morse and K.U. Ingard, \emph{Theoretical Acoustics},
McGraw-Hill 1968.

\bibitem{[PP96]} G. Propst and J. Pr\"uss, \emph{On wave equations with
boundary dissipation of memory type}, J. Integral Equations
Appl. \textbf{8} (1996), 99--123.

\end{thebibliography}
\end{document}